\documentclass[a4paper,11pt]{article}

\usepackage{amssymb}

\usepackage{amsmath}
\usepackage[all]{xy}

\usepackage[english,francais]{babel}

\newtheorem{e-proposition}[theorem]{Proposition}

\newtheorem{e-definition}[theorem]{Definition\rm}

\newtheorem{theoreme}{Th\'eor\`eme}[section]
\newtheorem{lemme}[theoreme]{Lemme}
\newtheorem{proposition}[theoreme]{Proposition}

\newtheorem{definition}[theoreme]{D\'efinition\rm}
\newtheorem{remarque}{\it Remarque}

\def\og{\leavevmode\raise.3ex\hbox{$\scriptscriptstyle\langle\!\langle$~}}
\def\fg{\leavevmode\raise.3ex\hbox{~$\!\scriptscriptstyle\,\rangle\!\rangle$}}

\begin{document}

\begin{center}
{\Large \bf Cohomologie du groupe lin\'{e}aire \`{a} coefficients dans les polyn\^{o}mes
de matrices}

\medskip

Antoine Touz\'{e}

\end{center}

\medskip 

\begin{abstract}
We compute bifunctors cohomology for matrix polynomials under conjugation and
detect candidates for universal classes in higher invariant theory.
\end{abstract}

\section{Introduction}
\label{}

Depuis une douzaine d'ann\'{e}es, la cohomologie des foncteurs a permis de
nombreux calculs et applications 
\cite{FS,FFSS,Chalupnik}, parmi lesquelles la d\'{e}montration
\cite{FS} de l'engendrement fini de l'alg\`{e}bre de cohomologie des
sch\'{e}mas en groupes finis. Dans \cite{FF}, Franjou et Friedlander
entament l'\'{e}tude de la cohomologie $H_\mathcal{P}^*(GL,B)$ des
bifoncteurs $B$ strictement polynomiaux. Cet article \'etend les
premiers calculs effectu\'{e}s dans \cite{FF}.

Nos r\'{e}sultats peuvent servir de base pour d'autres calculs et
laissent \'{e}galement entrevoir de nombreuses applications.
  Ainsi, la cohomologie d'un bifoncteur $B$ homog\`{e}ne de degr\'{e} 
$(b,b)$
calcule \cite[Th 1.5]{FF} la cohomologie rationnelle de $GL_n$ \`{a}
coefficients dans $B(k^n,k^n)$ pour $n\ge b$~:
$$H_\mathcal{P}^*(GL,B)\simeq H_{\mathrm{rat}}^*(GL_n,B(k^n,k^n))\;.$$

La cohomologie des bifoncteurs calcule \'{e}galement \cite[Th
8.2]{FF} la cohomologie stable des groupes discrets $GL_n(k)$ sur un
corps fini $k$. Comme il est not\'e dans \cite{EF} et \cite[section
8]{FF}, c'est l\`a un premier pas vers la K-th\'{e}orie de l'anneau
des nombres duaux et de $\mathbb{Z}/p^2\mathbb{Z}$.

\begin{theoreme}\label{thmprinc}
Soit $\mathbb{K}$ un corps de caract\'{e}ristique $p>0$ et $\mu =
(\mu_1,\dots,\mu_n)$ une partition de poids $d$. La s\'{e}rie de
Poincar\'{e} de
$H^*_\mathcal{P}(GL,S^{\mu_1(r)}gl\otimes\dots\otimes
S^{\mu_n(r)}gl)$ est \'{e}gale \`{a} la s\'{e}rie de Poincar\'{e}
des coinvarariants de $H^*_\mathcal{P}(GL,\otimes^{d(r)}gl)$ sous
l'action du sous-groupe
$\mathfrak{S}_{\mu_1}\times\dots\times\mathfrak{S}_{\mu_n}$ de
$\mathfrak{S}_d$.
\end{theoreme}

Ce th\'{e}or\`{e}me n'\'{e}tait auparavant connu que pour $d<p$
\cite[prop 4.1]{FF} et pour $d=p=2$ \cite[Th 5.1]{FF}. Il donne un
r\'{e}sultat ais\'{e}ment calculable \`{a} partir du
$\mathfrak{S}_d$-module $H^*_\mathcal{P}(GL,\otimes^{d(r)}gl)$ qui
est explicitement d\'{e}crit dans \cite[Th 1.8]{FF}. En particulier,
$H^*_\mathcal{P}(GL,S^{\mu(r)}gl)$ est nulle en degr\'{e}s impairs.

 Si l'on veut \'etendre les r\'esultats de \cite[\S
1]{FS}, van der Kallen explique dans \cite{vdK} l'int\'er\^et
d'obtenir des classes pour les coefficients en puissances divis\'ees
de $gl$.

\begin{proposition}
Soit $\mathbb{K}$ un corps de caract\'{e}ristique $p$. La s\'{e}rie de 
Poincar\'{e} de $H^*_\mathcal{P}(GL,\Gamma^{p(r)} gl)$ se d\'{e}duit
de celle de $H^*_\mathcal{P}(GL,S^{p(r)} gl)$ en ajoutant 
$(t^{2p-2}-1)\frac{1-t^{2p^{r+1}}}{1-t^{2p}}$.
\end{proposition}

En particulier, la s\'{e}rie de Poincar\'{e} de 
$H^*_\mathcal{P}(GL,\Gamma^{p(r)} gl)$ est nulle en degr\'{e}s impairs 
et elle a m\^{e}me
caract\'{e}ristique d'Euler-Poincar\'{e} que celle de $S^{p(r)}gl$. 
C'est en fait une cons\'{e}quence de ph\'{e}nom\`{e}nes de dualit\'{e}
  g\'{e}n\'{e}raux qui impliquent notamment le~:

\begin{theoreme}
Soit $F$ un foncteur strictement polynomial et $F^\sharp$ son dual de 
Kuhn. La
caract\'{e}ristique d'Euler-Poincar\'{e} de $H^*_\mathcal{P}(GL,F gl)$ 
est \'{e}gale
\`{a} celle de $H^*_\mathcal{P}(GL,F^\sharp gl)$.
\end{theoreme}

D'autres r\'{e}sultats donnant une description partielle de
$H^*_\mathcal{P}(GL,\Gamma^{n(r)} gl)$ peuvent \^{e}tre obtenus
\`{a} partir de notre connaissance de la cohomologie de
$S^{n(r)}gl$. Par exemple~:
\begin{proposition}
Le degr\'{e} maximal de la cohomologie de $\Gamma^{p^k(r)}gl$ est 
$p^k(2p^r-2)+2p^k-2$.
\end{proposition}

La d\'{e}monstration du th\'{e}or\`{e}me \ref{thmprinc} occupe la
deuxi\`{e}me partie de cet article. Elle repose sur un r\'{e}sultat
de th\'{e}orie des repr\'{e}sentations \cite{ABW}, un calcul d'
${\mathrm{Ext}}$ dans la cat\'{e}gorie $\mathcal{P}$
\cite{Chalupnik} et un argument de changement de corps. Pour plus de
clart\'{e}, le d\'{e}tail de cet argument est isol\'{e} dans la
troisi\`{e}me partie.

\section{Calcul de la s\'{e}rie de Poincar\'{e} de 
$H^*_\mathcal{P}(GL,S^{\mu(r)}gl)$}

Nous utilisons les notations suivantes.
Si $\mu=(\mu_1,\dots,\mu_n)$ est une partition et $E^*$ un foncteur 
exponentiel, $E^\mu$ d\'{e}signe le foncteur
$E^{\mu_1}\otimes\dots\otimes E^{\mu_n}$. Si $A$ et $B$ sont des 
foncteurs polynomiaux stricts, on note ${\mathcal{H}om}(A,B)$ le 
bifoncteur polynomial
strict ${\mathrm{Hom}}_\mathbb{K}(A(-_1),B(-_2))$. Les bifoncteurs de ce 
type sont dits s\'{e}parables et \cite[prop 2.2]{FF} donne un isomorphisme
naturel~: 
$H^*_\mathcal{P}(GL,{\mathcal{H}om}(A,B))\simeq{\mathrm{Ext}}^*_\mathcal{P}(A,B)$.

Le calcul se d\'{e}roule en deux \'{e}tapes. La premi\`{e}re
  repose sur une construction d'Akin, Buchsbaum et Weyman \cite[Th 
III.1.4 p.
244-245]{ABW} que nous reformulons ici en termes de bifoncteurs~:

\begin{theoreme}\label{theoremeABW}
Soit $k$ un entier. L'ordre lexicographique sur les partitions de poids 
$k$ induit une filtration de $S^k gl$ par des bifoncteurs strictement
polynomiaux~:
$$0\subseteq M_{(k)}\subseteq M_{(k-1,1)} \subseteq \dots \subseteq 
M_{(1,\dots,1)}=S^kgl$$
Le premier terme $M_{(k)}$ est isomorphe \`{a} 
${\mathcal{H}om}(\Lambda^k,\Lambda^k)$, et si $\lambda$ est une partition et
$\dot{\lambda}$ la suivante pour l'odre lexicographique on a une suite 
exacte courte~:
$$M_{\dot{\lambda}}\hookrightarrow M_{\lambda} \twoheadrightarrow 
{\mathcal{H}om}(W_{\lambda},S_{\lambda}) \;.$$
\end{theoreme}

La cohomologie du gradu\'{e} de cette filtration est bien connue 
\cite[Th 6.1]{Chalupnik}. En particulier, si $S_{\lambda/\lambda'}$
(resp. $W_{\lambda/\lambda'}$) d\'{e}signe le foncteur de Schur (resp. 
de Weyl) associ\'{e} au diagramme gauche $\lambda/\lambda'$,
la cohomologie du bifoncteur 
${\mathcal{H}om}({W_{\lambda/\lambda'}}^{(r)} 
,S_{\lambda/\lambda'}^{(r)} )$ est nulle en degr\'{e} impair.
En cons\'{e}quence, toutes les suites exactes courtes associ\'{e}es 
\`{a} la filtration de $S^{k(r)}gl$ scindent en cohomologie. La 
cohomologie de
${S}^{k(r)}gl$ est donc isomorphe \`{a} la cohomologie du gradu\'{e} de 
sa filtration.

Si $\mu$ est une partition, la cohomologie de $S^{\mu(r)}gl$ se calcule 
de mani\`{e}re similaire~: le produit
tensoriel des filtrations des $S^{\mu_i}gl$ donne une filtration de 
$S^\mu gl$ qui scinde en cohomologie. Pour exprimer le r\'{e}sultat
de fa\c{c}on concise, on note $(\lambda_1|\lambda_2|\dots|\lambda_n)$ le 
diagramme gauche tel que
$S_{(\lambda_1|\lambda_2|\dots|\lambda_n)}=S_{\lambda_1} 
\otimes\dots\otimes   S_{\lambda_n}$.

\begin{proposition}\label{prop-filtrationscinde}
Soit $\mu=(\mu_1,\dots,\mu_n)$ une partition de poids $d$. On a un 
isomorphisme d'espaces vectoriels gradu\'{e}s~:
$$ H^*_\mathcal{P}(GL,S^{\mu(r)}gl)\simeq 
\bigoplus_{|\lambda_1|=\mu_1\dots|\lambda_n|=\mu_n}
{\mathrm{Ext}}^*_\mathcal{P}
(W_{(\lambda_1|\lambda_2|\dots|\lambda_n)}^{(r)}
,S_{(\lambda_1|\lambda_2|\dots|\lambda_n)}^{(r)} ) $$
\end{proposition}

\begin{remarque}\label{rem-scinde}
Si la caract\'{e}ristique $p$ est grande, l'isomorphisme de la 
proposition peut s'obtenir directement en observant que la filtration 
scinde.
En effet, si $\lambda=(\lambda_1,\dots,\lambda_n)$ est une partition de 
poids $k<p$, l'application structurelle
$\Lambda^\lambda\twoheadrightarrow S_\lambda$ admet une
section $b_\lambda$. En notant $m:\otimes^k\twoheadrightarrow S^k$ le 
produit et $\Delta:\Lambda^\lambda\hookrightarrow\otimes^k$ le 
coproduit, on peut alors
d\'{e}finir $\phi_\lambda$ comme la compos\'{e}e~:
$${\mathcal{H}om}(W_\lambda,S_\lambda)
\xrightarrow[]{{\mathcal{H}om}(b_\lambda^\sharp,b_\lambda)}
{\mathcal{H}om}(\Lambda^\lambda,\Lambda^\lambda)\xrightarrow[]{{\mathcal{H}om}(\Delta^\sharp,\Delta)}
gl^{\otimes k}\xrightarrow[]{\frac{m}{\lambda_1!\dots\lambda_j!}} S^kgl\;.$$
L'examen de la construction d'Akin Buchsbaum et Weyman montre alors que 
$\oplus \phi_\lambda : \oplus {\mathcal{H}om}({W_{\lambda}} 
,S_{\lambda})\to S^k gl$ est un
isomorphisme.
\end{remarque}

Dans la suite, on note $\Delta^*$ le foncteur qui \`{a} un 
$\mathfrak{S}_d$-bimodule $M$ associe le $\mathfrak{S}_d$-module \`{a} 
gauche
obtenu par restriction de l'action de $\mathfrak{S}_d\times 
\mathfrak{S}_d^{\mathrm{op}}$
au sous-groupe 
$\mathfrak{S}_d=\{(\sigma,\sigma^{-1})\}\subset\mathfrak{S}_d\times\mathfrak{S}_d^{\mathrm{op}}$.
Si $f$ est un foncteur
des $\mathfrak{S}_d$-modules \`{a} gauche vers les $R$-modules et $M$ 
est un $\mathfrak{S}_d$-bimodule, on note $M f$
(resp $fM$) le $\mathfrak{S}_d$-module \`{a} gauche (resp. \`{a} droite) 
obtenu en appliquant $f$ \`{a} la structure de droite (resp. de
gauche).

La proposition \ref{prop-filtrationscinde} et \cite[Th.6.1]{Chalupnik} 
donnent une description de $H^*_\mathcal{P}(GL,S^{\mu(r)}gl)$
\`{a} partir du $\mathfrak{S}_d$-bimodule 
$H^*_\mathcal{P}(GL,\otimes^{d(r)}gl)$ d\'{e}crit dans \cite[Th.1.8]{FF} et
\cite[p.780]{Chalupnik}~:
$$H^*_\mathcal{P}(GL,S^{\mu(r)}gl)\simeq 
\bigoplus_{|\lambda_1|=\mu_1\dots|\lambda_n|=\mu_n}
\left(s_{(\lambda_1|\dots|\lambda_n)}\,H^*_\mathcal{P}(GL,\otimes^{d(r)}gl)\right)s_{(\lambda_1|\dots|\lambda_n)}\;.
$$
La deuxi\`{e}me \'{e}tape de la d\'{e}monstration consiste \`{a} 
simplifier cette expression. Pour cela, on peut tout d'abord supposer
qu'on travaille sur un corps de base $\mathbb{K}$ \'{e}gal \`{a} 
$\mathbb{F}_p$ d'apr\`{e}s \cite[Prop 3.1]{FF}. De plus, on remarque que
$H^*_\mathcal{P}(GL,\otimes^{d(r)}gl)$ est une somme directe de 
bimodules \'{e}l\'{e}mentaires dans le sens suivant~:

\begin{definition}
Soit $R$ un anneau. On appelle bimodule \'{e}l\'{e}mentaire sur $R$ un 
$R$-module de la forme~:
$R\mathfrak{S}_d/\mathfrak{S}_\gamma\otimes
R\mathfrak{S}_d$ avec $\mathfrak{S}_\gamma$ un sous-groupe de Young, 
muni de la structure de bimodule donn\'{e}e par l'action suivante sur la 
base~:
$$ \lambda.e_{[\tau]}\otimes e_{\sigma} .\mu = e_{[\lambda.\tau]}\otimes 
e_{\lambda.\sigma.\mu}\;. $$
\end{definition}

Lorsque $p$ est assez grand ($p>d$), on dispose d'une deuxi\`{e}me 
description de $H^*_\mathcal{P}(GL,S^{\mu(r)}gl)$. En effet la 
multiplication
$m:\otimes^d gl\to S^\mu gl$ induit un isomorphisme de $s^\mu \Delta^*\, 
H^*_\mathcal{P}(GL,\otimes^{d(r)}gl) $ vers 
$H^*_\mathcal{P}(GL,S^{\mu(r)}gl)$. Utilisons
les notations de la remarque \ref{rem-scinde}. D'apr\`{e}s \cite[Th 
6.1]{Chalupnik}, l'isomorphisme $\oplus
\phi_{\lambda_1}\otimes\dots\otimes\phi_{\lambda_n}$ induit une 
transformation naturelle~:
$$ \bigoplus_{|\lambda_1|=\mu_1\dots|\lambda_n|=\mu_n} 
(s_{(\lambda_1|\dots|\lambda_n)}\,-)\,s_{(\lambda_1|\dots|\lambda_n)}\to 
s^\mu \Delta^*$$
dont l'\'{e}valuation sur le $\mathfrak{S}_d$-bimodule 
${\mathrm{Ext}}^*_\mathcal{P}(\otimes^{d(1)},\otimes^{d(1)})$ est un 
isomorphisme.
Ce bimodule contient tous les bimodules \'{e}l\'{e}mentaires comme 
facteurs directs et on en d\'{e}duit donc une \'{e}galit\'{e} valable
pour tout bimodule \'{e}l\'{e}mentaire $M$ lorsque la 
caract\'{e}ristique $p$ est assez grande~:
$$\dim s^\mu \Delta^*\, M = \dim 
\bigoplus_{|\lambda_1|=\mu_1\dots|\lambda_n|=\mu_n} 
(s_{(\lambda_1|\dots|\lambda_n)}M)s_{(\lambda_1|\dots|\lambda_n)}\qquad(*)$$

Mais la dimension de $s^\mu \Delta^* \,M$ est ind\'{e}pendante de la 
caract\'{e}ristique, et d'apr\`{e}s la proposition 
\ref{prop-independancecaract},
il en va de m\^{e}me pour la dimension des $(s_\mu M)s_\lambda$. Ainsi, 
  l'\'{e}galit\'{e} $(*)$ est en fait valable en toute caract\'{e}ristique.

Le th\'{e}or\`{e}me d\'{e}coule maintenant du fait  que 
$H^*_\mathcal{P}(GL,\otimes^{d(r)}gl)$ est une somme directe de 
bimodules \'{e}l\'{e}mentaires.

\section{Changement de base pour les bimodules \'{e}l\'{e}mentaires }

\begin{lemme}
Soit $M$ un bimodule \'{e}l\'{e}mentaire sur $\mathbb{F}_p$ et $\mu$ et 
$\lambda$ deux diagrammes gauches. Les applications naturelles~:
$(s_\mu M)_{\mathfrak{S}_\lambda}\to s_\mu (M_{\mathfrak{S}_\lambda}) $ 
et $ s_\mu (M^{{\mathrm{alt}}\mathfrak{S}_\lambda}) \to  (s_\mu 
M)^{{\mathrm{alt}}\mathfrak{S}_\lambda}$ sont des isomorphismes.
\end{lemme}

\noindent\textit{D\'{e}monstration.}
Nous d\'{e}montrons le premier isomorphisme. On a un diagramme commutatif~:
$$\xymatrix{
(s_\mu M)_{\mathfrak{S}_\lambda}\ar@{->}[r]^{(1)}
& s_\mu (M_{\mathfrak{S}_\lambda})\\
(\,^{{\mathrm{alt}}\mathfrak{S}_\mu} 
M)_{\mathfrak{S}_\lambda}\ar@{->}[r]^{(2)}\ar@{->}[u]
& \,^{{\mathrm{alt}}\mathfrak{S}_\mu} 
(M_{\mathfrak{S}_\lambda})\ar@{->>}[u]\\
}$$
La fl\`{e}che $(2)$ est un isomorphisme car les actions sont 
d\'{e}finies sur une base de  $M$. Par cons\'{e}quent, $(1)$ est surjective.
On peut conclure car la source et le but de $(1)$ ont la m\^{e}me 
dimension d'apr\`{e}s \cite[Th 6.1]{Chalupnik}.

\begin{definition}
Soit $M$ un $\mathfrak{S}_d$-bimodule et $\lambda$ et $\mu$ des 
diagrammes gauches. On d\'{e}finit $s_{\mu} M s_{\lambda} $ comme 
l'image de l'application
compos\'{e}e~:
$$^{\mathfrak{S}_\mu{\mathrm{alt}}} 
(M^{{\mathrm{alt}}\mathfrak{S}_\lambda})\hookrightarrow M \twoheadrightarrow
(\,_{\mathfrak{S}_{\widetilde{\mu}}}
M)_{\mathfrak{S}_{\widetilde{\lambda}}}\;. $$
\end{definition}

\begin{lemme}\label{lmineg}
Soit $M$ un bimodule \'{e}l\'{e}mentaire sur  $\mathbb{Z}$ et $\lambda$ et
$\mu$ des diagrammes gauches. On a~:
$$\dim \left(s_\mu \left(\mathbb{F}_p\otimes M\right)\right)s_\lambda 
\le \mathrm{rank}\, s_{\mu} M s_\lambda$$
De plus, l'in\'{e}galit\'{e} est en fait une \'{e}galit\'{e} si $p$ est 
assez grand.
\end{lemme}

\noindent\textit{D\'{e}monstration.}
Notons $M_{\mathbb{F}_p}$ le bimodule \'{e}l\'{e}mentaire 
$M\otimes\mathbb{F}_p$.
On a un diagramme commutatif~:
$$\xymatrix{
^{\mathfrak{S}_\mu{\mathrm{alt}}} 
(M_{\mathbb{F}_p}^{{\mathrm{alt}}\mathfrak{S}_\lambda})\ar@{^{(}->}[d]\ar@{->>}[r]
& s_\mu(M_{\mathbb{F}_p}^{{\mathrm{alt}}\mathfrak{S}_\lambda}) 
\ar@{^{(}->}[d]\ar@{->}[r]^{\simeq}
& (s_\mu M_{\mathbb{F}_p})^{{\mathrm{alt}}\mathfrak{S}_\lambda} 
\ar@{^{(}->}[d]\ar@{->>}[r]
& (s_\mu M_{\mathbb{F}_p})_{\mathfrak{S}_{\widetilde{\lambda}}} 
\ar@{->}[d]^{{i}_{\mathfrak{S}_{\widetilde{\lambda}}}}\\
M_{\mathbb{F}_p}\ar@{->>}[r]
& _{\mathfrak{S}_{\widetilde{\mu}}}{M_{\mathbb{F}_p}} \ar@{->}[r]^{=}
& _{\mathfrak{S}_{\widetilde{\mu}}}{M_{\mathbb{F}_p}}\ar@{->>}[r]
& (\,_{\mathfrak{S}_{\widetilde{\mu}}} 
M_{\mathbb{F}_p})_{\mathfrak{S}_{\widetilde{\lambda}}}\\
}$$
On v\'{e}rifie que ${i}_{\mathfrak{S}_{\widetilde{\lambda}}}$ est 
injective. Une chasse dans le diagramme montre alors que
$s_{\mu} M_{\mathbb{F}_p} s_\lambda =\left(s_\mu M_{\mathbb{F}_p} 
\right)s_\lambda $.

Regardons maintenant le diagramme de $\mathbb{Z}$-modules dont les 
fl\`{e}ches verticales sont des surjections
  induites par la r\'{e}duction modulo $p$~:
$M\twoheadrightarrow M_{\mathbb{F}_p}$~:
$$\xymatrix{
^{\mathfrak{S}_\mu{\mathrm{alt}}} 
(M^{{\mathrm{alt}}\mathfrak{S}_\lambda})\ar@{->>}[d]\ar@{^{(}->}[r]
& M \ar@{->>}[d]\ar@{->>}[r] & 
(\,_{\mathfrak{S}_{\tilde{\mu}}}M)_{\mathfrak{S}_{\tilde{\lambda}}} 
\ar@{->>}[d]\\
^{\mathfrak{S}_\mu{\mathrm{alt}}} 
(M_{\mathbb{F}_p}^{{\mathrm{alt}}\mathfrak{S}_\lambda})\ar@{^{(}->}[r]& 
M_{\mathbb{F}_p} \ar@{->>}[r]
& 
(\,_{\mathfrak{S}_{\tilde{\mu}}}M_{\mathbb{F}_p})_{\mathfrak{S}_{\tilde{\lambda}}} 
\\
}$$
Puisque 
$(\,_{\mathfrak{S}_{\widetilde{\mu}}}M)_{\mathfrak{S}_{\widetilde{\lambda}}}$ 
est un $\mathbb{Z}$-module libre, il en va de m\^{e}me pour
$s_\mu Ms_\lambda$ et on en d\'{e}duit le r\'{e}sultat.

\begin{lemme}\label{lm-SFB}
Soient $\lambda$ et $\mu$ des diagrammes gauches et $k$ un entier.
La dimension de ${\mathrm{Hom}}_{\mathcal{P}} (W_\mu\circ kI,S_\lambda)$ 
ne d\'{e}pend pas du corps.
\end{lemme}

\noindent\textit{D\'{e}monstration.}
Gr\^{a}ce \`{a} \cite[Th.II.2.16 et II.3.16]{ABW}, les foncteurs de 
Schur $S_\lambda$ et de Weyl $W_\mu$
sont en fait des foncteurs polynomiaux stricts sur $\mathbb{Z}$ au sens 
de \cite[section 2]{SFB}.
En copiant la d\'{e}monstration de \cite[Prop.2.8]{SFB}, on montre que 
le foncteur ${S_\lambda}\otimes_\mathbb{Z}{\mathbb{F}_p}$
(resp. ${W_\mu}\otimes_\mathbb{Z}{\mathbb{F}_p}$)
obtenu par changement de base n'est autre que le foncteur de Schur 
(resp. Weyl) d\'{e}fini sur ${\mathbb{F}_p}$.
On utilise alors \cite[Prop.2.6]{SFB} pour conclure.

\begin{proposition}\label{prop-independancecaract}
Soient $\lambda$, $\mu$ des diagrammes gauches de poids  $d$. Pour tout 
bimodule \'{e}l\'{e}mentaire
   $M=\mathbb{F}_p\mathfrak{S}_n/\mathfrak{S}_\gamma\otimes 
\mathbb{F}_p\mathfrak{S}_n$ la dimension de
$(s_\mu M)s_\lambda$ ne d\'{e}pend pas de $p$.
\end{proposition}

\noindent\textit{D\'{e}monstration.}
Les d\'{e}monstrations de \cite[Th.4.4 et Th 6.1]{Chalupnik} montrent 
l'\'{e}galit\'{e} suivante, valable pour tout entier $k$ et tout nombre
premier $p$~:
$$ \dim {\mathrm{Hom}}_{\mathcal{P},\mathbb{F}_p}(W_\mu\circ kI,S_\lambda) =
\dim 
\left(s_\mu\,{\mathrm{Hom}}_{\mathcal{P},\mathbb{F}_p}((kI)^{\otimes 
d},\otimes^d)\right)s_\lambda\qquad(*)$$

Soit $p$ et $q$ des nombres premiers, avec $q$ suffisamment grand pour 
que l'in\'{e}galit\'{e} du lemme \ref{lmineg}
soit une \'{e}galit\'{e} pour tout $\mathfrak{S}_d$-bimodule 
\'{e}l\'{e}mentaire. (Ceci est possible car il n'y en a qu'un nombre 
fini \`{a} isomorphisme pr\`{e}s).
Alors pour tout bimodule \'{e}l\'{e}mentaire $M$ sur $\mathbb{Z}$ on a 
l'in\'{e}galit\'{e}~:
$$\dim (s_\mu \left(M\otimes\mathbb{F}_p\right))s_\lambda \le \dim 
(s_\mu (M\otimes\mathbb{F}_q))s_\lambda \;.$$
Mais ${\mathrm{Hom}}_{\mathcal{P},\mathbb{F}_p}((dI)^{\otimes 
d},\otimes^d)$ contient tous les bimodules \'{e}l\'{e}mentaires comme 
facteurs directs.
En cons\'{e}quence, $(*)$ et le lemme \ref{lm-SFB} montrent que 
l'in\'{e}galit\'{e} pr\'{e}c\'{e}dente est en fait une \'{e}galit\'{e}, 
ce qui ach\`{e}ve la
d\'{e}monstration.

\end{document}